\documentclass{article}

\usepackage{amsmath, amsthm, amsfonts}
\usepackage{color}

\usepackage{hyperref}

\usepackage{longtable}

\usepackage{csquotes}

\newtheorem{thm}{Theorem}[section]
\newtheorem{cor}[thm]{Corollary}
\newtheorem{lem}[thm]{Lemma}
\newtheorem{theorem}[thm]{Theorem}

\newtheorem{defn}[thm]{Definition}

\usepackage{mdframed}
\theoremstyle{definition}
\newmdtheoremenv{boxProb}{Problem}
\newmdtheoremenv{boxDef}{Definition}
\newmdtheoremenv{boxCor}{Corollary}
\newmdtheoremenv{boxThm}{Theorem}
\newmdtheoremenv{compjob}{Computational Job}
\newmdtheoremenv{reqi}{Requirement}

\usepackage{graphicx}
\usepackage{subfigure}          

\usepackage{placeins}

\usepackage{tabularx}

\usepackage{lmodern,textcomp}

\usepackage{algorithm}
\usepackage{algpseudocode}

\usepackage{enumitem}

\usepackage{xspace} 
\newcommand\largeparbreak{\par\bigskip}

\newcommand{\inv}{^{-1}\xspace}

\newcommand{\bxi}{\boldsymbol{\xi}\xspace}

\renewcommand{\t}{^\textsf{T}\xspace}

\newcommand{\dom}{\operatorname{dom}\xspace}

\newcommand{\away}[1]{}

\newcommand{\R}{\mathbb{R}\xspace}

\newcommand{\N}{\mathbb{N}\xspace}

\newcommand{\cB}{\mathcal{B}\xspace}

\newcommand{\cD}{\mathcal{D}\xspace}

\newcommand{\cJ}{\mathcal{J}\xspace}

\newcommand{\cS}{\mathcal{S}\xspace}

\newcommand{\cO}{\mathcal{O}\xspace}

\newcommand{\bQ}{\mathbf{Q}\xspace}

\newcommand{\bx}{\mathbf{x}\xspace}
\newcommand{\bxL}{\mathbf{x}_{\textsf{L}}\xspace}
\newcommand{\bxR}{\mathbf{x}_{\textsf{R}}\xspace}

\newcommand{\bc}{\mathbf{c}\xspace}

\newcommand{\hbx}{\hat{\mathbf{x}}\xspace}

\newcommand{\be}{\mathbf{1}\xspace}
\newcommand{\bei}[1]{{\mathbf{e}}\xspace}

\newcommand{\tol}{{\textsf{tol}}\xspace}

\newcommand{\NCQPTR}{\textsf{NCQPBTR}}

\newcommand{\GammLR}{{\Gamma_{\textsf{LR}}}\xspace}
\newcommand{\GammDelta}{{\Gamma_{\Delta}}\xspace}
\newcommand{\GammD}{{\hat{\Gamma}}\xspace}

\usepackage{todonotes}

\title{Weakly polynomial efficient minimization of a non-convex quadratic function with logarithmic barriers in a trust-region}

\author{Martin Neuenhofen}

\begin{document}

\maketitle

\begin{abstract}
We introduce a particular optimization problem that minimizes the sum of a non-convex quadratic function and logarithmic barrier-functions in a $\ell_\infty$-trust-region (i.e. cube). Our paper covers three topics.

We explain the relevance of the considered problem. We lay out how solutions of this problem can be used as efficient step-directions in solution methods for nonlinear programming.

We present a theoretical algorithm for solving the problem. We show that this algorithm has weak polynomial time-complexity.

A practical method is under development. In the outlook we discuss how the given method can be accelerated for better practical performance. We also lay out where the difficulties live when trying to formulate an accelerated primal-dual variant.
\end{abstract}


\section{Introduction}
\subsection{Problem statement}

We consider the minimization of the following program that, as we motivate later, we call \textit{non-convex quadratic program with barriers and trust-region} (\NCQPTR):
\begin{equation}
	\tag{\NCQPTR}\label{eqn:NCQPTR}
	\begin{aligned}
		&\operatornamewithlimits{min}_{\bx \in \Omega}& \Phi(\bx) := &\quad \frac{1}{2} \cdot \bx\t \cdot \bQ \cdot \bx + \bc\t \cdot \bx \\
		& & &\quad - \tau_F \cdot \be\t \cdot \Big(\,\log(\bx-\bxL) + \log(\bxR-\bx)\,\Big) \\
		& & &\quad - \pi_F \cdot \be\t \cdot \Big(\,\log(\Delta\cdot\be+\bx) + \log(\Delta\cdot\be-\bx)\,\Big)\,.
	\end{aligned}
\end{equation}
$\be \in \R^n$ is the vector of all ones. In this problem, the symmetric matrix $\bQ \in \R^{n \times n}$, vectors $\bc,\bxL,\bxR \in \R^n$, and positive scalar parameters $\tau_F \geq \pi_F>0$, $\Delta > 0$ are given. The solution domain is defined as follows:
\begin{align*}
	&\text{Box:} \quad& \cB &:= \lbrace \,\bxi \in \R^n \, \vert \, \bxL < \bxi < \bxR \, \rbrace\,;\\
	&\ell_\infty\text{-Sphere:} \quad& \cS &:= \lbrace \,\bxi \in \R^n \, \vert \, -\Delta\cdot\be < \bxi < \Delta\cdot\be \, \rbrace\,;\\
	&\text{Domain:} \quad&\Omega &:= \cB \cap \cS\,.
\end{align*}
We notice that $\overline{\Omega}$ is a box. We write the scalar parameter $\delta$ to denote the shortest side-length of this box.

We state two conditions. These must be satisfied in order for our proposed algorithm to be able of finding a solution in weakly polynomial time-complexity. The conditions are:
\begin{enumerate}
	\item It must be $\delta > 0$, i.e. $\Omega$ must be non-empty.
	\item The subsequently defined function $\psi$ must be convex over $\Omega$.
\end{enumerate}
The function $\psi\,:\,\Omega\rightarrow\R$ is defined below.
\begin{align*}
	\psi(\bx) := \frac{1}{2} \cdot \bx\t \cdot \bQ \cdot \bx + \bc\t \cdot \bx - \frac{\tau_F}{2} \cdot \be\t \cdot \Big(\,\log(\bx-\bxL) + \log(\bxR-\bx)\,\Big)\,.
\end{align*}

Remark: We stated the first condition deliberately in terms of $\delta$ because the complexity of our algprithm grows with $1+\log_2\big(1-\log_2(\delta)\big)$. Hence, $\delta$ must be bounded below. Otherwise, (very) small values of $\delta$ could blow up the complexity of our method.

\subsection{Interpretation of the problem}
Problem \eqref{eqn:NCQPTR} can be interpreted as an ($\bxL,\bxR$)-box-constrained non-convex quadratic program in a $\ell_\infty$-trust-region of radius $\Delta>0$\,. Since $\overline{\Omega}$ is actually just another box, \eqref{eqn:NCQPTR} still resembles a box-constrained non-convex quadratic program when $\tau_F,\pi_F \searrow 0$.

Unfortunately, in general it is NP-hard to solve a box-constrained non-convex quadratic program \cite{478584,MR1770270}. But, there are two significant differences between \eqref{eqn:NCQPTR} and the general box-constrained non-convex quadratic program. These differences allow us to solve \eqref{eqn:NCQPTR} in an efficient way. We discuss these differences below.

\largeparbreak

First, it is an important detail that in our numerical setting the value $\tau_F>0$ is a finite positive value. I.e., there is no convergence of $\tau_F$ to $+0$. Hence, the logarithmic barriers with $\bxL,\bxR$ introduce a crucial amount of convexity to $\psi$.

Second, while the box-constrained non-convex quadratic program considers a general matrix $\bQ$, we made the particular requirement that $\psi$ must be convex over $\Omega$. This imposes a restriction on $\bQ$ for the problem that we consider here. Clearly, this restriction is parametric in $\tau_F$, as sufficiently large values of $\tau_F$ can always force convexity of $\psi$, regardless of $\bQ$.

\subsection{Motivation of the problem}
Below, we motivate the merits of using solutions to problems \eqref{eqn:NCQPTR} as a tool for constructing an iterative scheme that solves non-convex nonlinear programming problems (NLP). Our motivation must be superficial and simplistic, as it actually goes deeply into the convergence mechanisms of numerical methods for non-convex nonlinear programming.

In our motivation, we are superficial and simplifying when it comes to discussing the following: approaches for building an iteration scheme, the strengths and disadvantages of a respective iteration scheme, and judgments on practical suitability of an iteration scheme. We also do not discuss equality constraints.

\largeparbreak

As a simplification, let us agree that solvers for large-scale nonlinear programming can be divided into two classes: successive quadratic programming (SQP) and interior-point methods (IPM). The first class is characterized by solving a quadratic program to determine the step-direction (cf., e.g., in \cite{WORHP}), while the second class determines the step-direction by solving some sort of barrier-regularized Karush-Kuhn-Tucker (KKT) equations with Newton's method (cf., e.g., in \cite{ForsgrenGill,ChenGoldfarb}).

To promote global convergence, both methods require an inertia correction. (For the scope of this discussion we do not care whether this inertia correction is interpreted as Quasi-Newton line-search or trust-region approach.) In the general case, when the iteration of either SQP or IPM approaches a local minimum, the inertia in the Hessian matrix of SQP and IPM behave entirely different:

\paragraph{In SQP}, close to (and even in) the limit point of the iteration, the quadratic program for the step-direction may be non-convex. Hence, in order to compute a step-direction in SQP, there are two options: 

As one option, one solves the non-convex quadratic program as is. However, superficially said, this is impractical because the solution of non-convex quadratic programming is NP-hard (even when only solving it for a local minimizer) \cite{TrustRegionMethods}. Thus, each iteration could be prohibitively expensive, which must be avoided.

As an alternative, the quadratic program can be made convex (using a shift or a positive definite Hessian approximation), and can then be solved efficiently. This approach is followed, e.g., in \cite{SNOPT,WORHP}. But, the convexization of the quadratic program impedes the accuracy in which the nonlinear program is locally approximated. Hence, as a serious drawback, in the general case these methods do not converge locally fast.

\paragraph{In IPM}, in contrast, in the vicinity of the limit point of the iteration, the Newton-system for the barrier-regularized KKT equations has the correct inertia. This is easily explained from the fact that IPM replace the inequality constraints with barriers. Thus, the limit point of the iteration minimizes an unconstrained function. Since locally minimal, that unconstrained function must be locally convex, and hence the matrix has the correct inertia automatically.

But, there is also an issue with IPM. Namely, the Newton iteration suffers slow convergence if the iterates do not live close to a curve in space that is referred to as central path; see in \cite{IniGuessIPM} for a discussion. This is a serious drawback, as it constitutes a non-trivial restriction on the choice of an initial guess for starting the iteration.

\paragraph{New approach}
Our method can improve on IPM because it allows us to milden the restriction on the initial guess, that IPM imposes. This can work as follows: We use an IPM principle, i.e. we solve a sequence of barrier-problems for decreasing values of a barrier-parameter $\tau_{\text{IPM}}>0$.

For the step-direction, we use a solution of \eqref{eqn:NCQPTR}, where $\tau_F$ is chosen as the current value of the barrier-parameter $\tau_{\text{IPM}}$ and where the NLP objective is approximated with a quadratic function through $\bQ$ and $\bc$. Globalization of the IPM iteration is achieved with a trust-region of radius $\Delta$. This same value of $\Delta$ shall be used in \eqref{eqn:NCQPTR} when computing the step-direction.

Building on the strengths of IPM, in the limit the function $\psi$ will be convex automatically, hence the matrix $\bQ$ can be chosen as the (possibly non-convex) Hessian of the NLP objective. Clearly, the iteration converges locally rapid, because in the limit $\psi$ can be chosen as a second-order approximation to the objective of the barrier-subproblem, and $\pi_E$ can be chosen orders of magnitude smaller.

\section{Preliminaries}
Because for our method we present a rigorous analysis, we make use of a couple of non-trivial results from the literature. These results deal with properties of convex functions and they deal with the iteration complexity of Newton-type methods applied to minimizing these functions. All these literature results are given here in one place for the reader's convenience.

\subsection{Notation and definitions}

We use the symbol $q$ for a quadratic function and the symbol $\Gamma$ for a logarithmic barrier-function. Both map from $\R^n$ into $\R$. Logarithmic barrier-functions are functions of the form
\begin{align*}
	\Gamma(\bx) := -\sum_{k =1}^d \sum_{j \in \cJ_k} r_{j,k} \cdot \log(c_{j,k} + s_{j,k} \cdot x_j)\,,
\end{align*}
where $\cJ_k \subset \lbrace 1,...,n\rbrace$ for $k=1,..,d$, and $d \in \N$. There are numbers $c_{j,k} \in \R$, signs $s_{j,k} \in \lbrace -1,1\rbrace$, and scalars $r_{j,k} \in [1,\infty)$. The value $x_j$ is the $j$th component of $\bx$. Below, we give some particular instances for such functions $q,\Gamma$.

\largeparbreak

We define the quadratic function of our optimization problem for the given user-inputs $\bQ,\bc$:
\begin{align}
	\hat{q}(\bx) := \frac{1}{2} \cdot \bx\t \cdot \bQ \cdot \bx\t + \bc\t \cdot \bx\label{eqn:hatq}
\end{align}
We define the barrier-function for the box $\cB$
\begin{align*}
	\GammLR(\bx) := - \be\t \cdot \Big(\,\log(\bx-\bxL) + \log(\bxR-\bx)\,\Big)\,,
\end{align*}
and for the trust-region $\cS$
\begin{align*}
	\GammDelta(\bx) := - \be\t \cdot \Big(\,\log(\Delta\cdot\be+\bx) + \log(\Delta\cdot\be-\bx)\,\Big)\,,
\end{align*}
and for the solution domain $\Omega$
\begin{align}
	\GammD(\bx) := \GammLR(\bx) + \GammDelta(\bx)\,.\label{eqn:hatGamma}
\end{align}

An important framework, that our entire analysis is based on, is that of self-concordance \cite{Boyd,Nesterov}.
\begin{defn}[Self-concordance]
	Let $g : \dom(g) \subset \R^d \rightarrow \R$ be three times differentiable, where $d=1$. Let $\cD \subset \dom(g)$ be convex, bounded, and open. Then $g$ is called {self-concordant} over $\cD$ if it satisfies
	\begin{align*}
		|g'''(x)| \leq 2 \cdot \big(\,g''(x)\,\big)^{\frac{3}{2}}\quad \forall x \in \cD\,.
	\end{align*}
	The definition is generalized to the case $d\in\N$ by requiring that the above condition holds for every projection of $g$ and $\cD$ onto a straight line through $\cD$.
\end{defn}
\noindent
It is well understood that logarithmic barrier-functions and convex quadratic functions are self-concordant. The right-hand side in the condition necessitates convexity of $g$ in $\cD$.

We make a further definition, that serves as a measure for the amount of constraints that a logarithmic barrier-function reflects.
\begin{defn}[Absolute value of a barrier-function]
	Consider a logarithmic barrier-function $\Gamma$. We define the absolute value of $\Gamma$ as
	\begin{align*}
		|\Gamma|:= \sum_{k =1}^d \sum_{j \in \cJ_k} r_{j,k}\,.
	\end{align*}
\end{defn}
\noindent
The following hold:
\begin{align*}
	|\GammLR| = 2 \cdot n\,,\quad\quad
	|\GammDelta| = 2 \cdot n\,,\quad\quad
	|\GammD| = 4 \cdot n\,.
\end{align*}

\subsection{Tools for minimizing self-concordant functions}
The framework of self-concordance originated from work of Nesterov \cite{Nesterov}. The conduction of analysis for self-concordant functions in the literature is caused by their desirable features. In particular, as the major desirable feature of self-concordant functions, they can be minimized efficiently with Newton's method.

Our proposed algorithm is based entirely on this capability of Newton's method. We review two crucial results in this regard, that will be needed in our analysis.

\subsubsection{Newton's method for minimization}
The first result considers the situation when we apply a damped version of Newton's method to minimize a self-concordant function $g$ from an initial guess $\bx_0 \in \dom(g)$. The theorem and algorithm below are summarized from \cite[Section~9]{Boyd}.

\begin{theorem}\label{thm:Newton_Minimization}
	Let $g : \dom(g) \rightarrow \R$ be strictly convex and $\dom(g)$ bounded. Further, let $g$ be self-concordant over $\dom(g)$. Let $\bx_0 \in \dom(g)$. Choose $\varepsilon>0$.
	
	Then, after $k \in \N$ iterations of Algorithm~\ref{algo:DampedNewton}, the iterate $\bx_k$ satisfies
	\begin{align*}
		g(\bx_k)-\min\lbrace g \rbrace \leq \varepsilon\,,
	\end{align*}
	where $k$ obeys to the bound
	\begin{align*}
	k \leq 375 \cdot \big(\,g(\bx_0) - \min\lbrace g\rbrace\,\big) + \log_2\big(1-\log_2(\varepsilon)\big)\,.
	\end{align*}
	Furthermore, if $\bx_k$ for an arbitrary $k \in \N_0$, satisfies
	$$ 	g(\bx_k)-\min\lbrace g\rbrace \leq 0.25\,,	$$
	then the Newton iteration terminates in the next $\lceil\log_2(1-\log_2(\varepsilon))\rceil$ iterations.
\end{theorem}

\begin{algorithm}
	\caption{Damped Newton method}
	\label{algo:DampedNewton}
	\begin{algorithmic}[1]
		\Procedure{DampedNewton}{$g,\bx_0,\varepsilon$}
		\For{$k=0,1,2,\dots$}
			\State $\Delta\bx_k := -\nabla^2 g(\bx_k)\inv \cdot \nabla g(\bx_k)$
			\State $\lambda_k^2 := -\Delta\bx\t_k \cdot \nabla g(\bx_k)$
			\If{$\lambda^2 / 2 \leq \varepsilon$}
				\State \Return $\bx_k$
			\EndIf
			\State $t:=1$
			\While{$g(\bx_k + t \cdot \Delta\bx_k) > g(\bx_k) - 0.1 \cdot t \cdot \lambda_k^2$}
				\State $t := 0.8 \cdot t$
			\EndWhile
			\State $\bx_{k+1} := \bx_k + t \cdot \Delta\bx_k$
		\EndFor
		\EndProcedure
	\end{algorithmic}
\end{algorithm}

The above result is useful when we have an initial guess $\bx_0$, whose function value is close to $\min\lbrace g\rbrace$. But sometimes we do not have such a point. In this situation, it may be more helpful if we can bound the number $k$ of Newton iterations by the value of $\|\nabla g(\bx_0)\|_2$. This is also meaningful, since, at a qualitatively good initial guess for the minimizer, the gradient should be close to zero.

In \cite{Boyd} the authors show that the \emph{Newton decrement} $\lambda$, that is computed in Algorithm~\ref{algo:DampedNewton}, serves as a bound for the optimality gap. In particular, it holds
\begin{align}
	g(\bx) - \min\lbrace g \rbrace \leq \lambda(\bx)^2 /2\label{eqn:BoundOptimalityGap}
\end{align}
where the Newton decrement $\lambda$ is defined as
\begin{align*}
	\lambda(\bx) := \sqrt{\nabla g(\bx)\t \cdot \nabla^2 g(\bx)\inv \cdot \nabla g(\bx)}
\end{align*}
From the definition we find the bound
\begin{align}
	\lambda(\bx) &\leq \sqrt{\frac{1}{\lambda_{\text{min}}(\nabla^2 g(\bx))}} \cdot \|\nabla g(\bx)\|_2\,.\label{eqn:BoundNewtonDecrement}
\end{align}

\subsubsection{Newton iterations for primal path-following}
Now, in this subsection, we consider a different result, where this time Newton's method is applied to minimize a parametric function.

The result can be used when we want to minimize a parametric self-concordant function with barrier terms of the following form:
\begin{align}
	\phi_\tau(\bx) := \frac{1}{\tau} \cdot f(\bx) + \Gamma(\bx)\,,\quad\quad \tau >0\,,
\end{align}
where $\phi_{\tau}$ is a self-concordant function for all relevant values of $\tau$ (, as defined below in detail).

Often, we have a minimizer $\bx_0$ for $\phi_{\tau_0}$, where $\tau_0$ is a moderate value.
But, we actually search for a minimizer $\bx_E$ of $\phi_{\tau_E}$, where $\tau_E$ is significantly smaller than $\tau_0$. The theorem and algorithm below show an efficient way to compute $\bx_E$ from $\bx_0$. The theorem and the algorithm are summarized from \cite[Section~11]{Boyd}.
\begin{theorem}\label{thm:NewtonPathFollow}
	Consider $\phi_\tau$ from above. Let $\dom(\Gamma)$ be bounded, and let $\phi_\tau$ be strictly convex and self-concordant in $\dom(\phi) := \dom(f) \cap \dom(\Gamma)$ for all $\tau \in [\tau_E,\tau_0] \subset (0,\infty)$. Let $\bx_0 \in \dom(\phi)$ satisfy
	\begin{align*}
		\phi_{\tau_0}(\bx_0) - \min\lbrace \phi_{\tau_0} \rbrace \leq 1/4\,.
	\end{align*}
	Choose $\varepsilon >0$ and $\tau_E \in (0,\tau_0)$.
	\largeparbreak
	Then, after $j \in \N$ iterations of Algorithm~\ref{algo:PathFollowNewton}, a solution $\bx_E$ is returned that satisfies
	\begin{align*}
		\phi_{\tau_E}(\bx_E) - \min\lbrace \phi_{\tau_E} \rbrace \leq \varepsilon\,.
	\end{align*}
	The number of outer iterations $j$ is bounded by
	\begin{align*}
		j \leq \left\lceil\frac{\log\big(\,|\Gamma| \cdot \tau_0/\tau_E\,\big)}{\log(1+1/\sqrt{m})}\right\rceil \in \cO\Big(\,\sqrt{|\Gamma|} \cdot \big(\log|\Gamma|+\log(\tau_0/\tau_E)\big)\,\Big)\,.
	\end{align*}
	The number of each inner iterations of Algorithm~\ref{algo:DampedNewton} during the for-loop in line~6 is bounded by
	\begin{align*}
		k \leq 380\,.
	\end{align*}
	Only the final call of Algorithm~\ref{algo:DampedNewton} in line~12 requires a number $k \in \N$ of iterations, that is not bounded by a constant. Namely, instead it is bounded as
	\begin{align*}
		k \leq \left\lceil\log_2\big(1-\log_2(1/\varepsilon)\big)\right\rceil\,.
	\end{align*}
	
	Hence, in summary, the total number of linear equation systems, that must be solved in Algorithm~\ref{algo:PathFollowNewton}, lives in
	\begin{align*}
		\cO\Big(\, \sqrt{|\Gamma|} \cdot \big(\log|\Gamma|+\log(\tau_0/\tau_E)\big) + \log\big(1-\log(\varepsilon)\big) \,\Big)\,.
	\end{align*}
\end{theorem}

\begin{algorithm}
	\caption{Path-following Newton method}
	\label{algo:PathFollowNewton}
	\begin{algorithmic}[1]
		\Procedure{PathFollowNewton}{$f,\Gamma,\bx^{(0)},\tau_0,\tau_E,\varepsilon$}
		\State $\sigma := \frac{1}{1+1/\sqrt{|\Gamma|}}$\,,\quad$\tau := \tau_0$
		\For{$j=0,1,2,\dots$}
			\State $\tau := \max\lbrace\,\tau_E\,,\,\sigma \cdot \tau\,\rbrace$
			\State $g(\cdot) := 1/\tau \cdot f(\cdot)+\Gamma(\cdot)$\,,\quad$\bx_0 := \bx^{(j)}$
			\State $\bx_k:= $\textsc{DampedNewton}($g,\bx_0,1/4$)
			\State $\bx^{(j+1)} := \bx_k$
			\If{$\tau == \tau_E$}
				\State \textbf{break for-loop}
			\EndIf
		\EndFor
		\State\Return $\bx_E := $\textsc{DampedNewton}($g,\bx^{(j+1)},\varepsilon$)
		\EndProcedure
	\end{algorithmic}
\end{algorithm}

\subsection{Tools for composing self-concordant functions}
In the former subsection we have seen that for self-concordant functions there exist simple methods with strong theoretical results. These results provide upper bounds for the iteration complexity of Newton's method for minimization.

In this subsection we give results that can be used to construct a self-concordant function from a convex function.
\begin{lem}\label{lem:self_from_conv}
	Consider the function $q(\bx)+f(\bx)$, where $q$ is quadratic and $f$ is a self-concordant function with $\dom(f) \subset \R^d$, and $d \in \N$. Let $q+f$ be convex over $\cD \subset \dom(f)$.
	
	Then, the function
	\begin{align*}
		g(x) := 8 \cdot \big(\,q(\bx) + 2 \cdot f(\bx)\,\big)
	\end{align*}
	is self-concordant over $\cD$.
\end{lem}
\noindent
\underline{Proof:}
Meaning all the following relations on a fixed arbitrary straight line through $\cD$, we find from the self-concordance of $f$:
\begin{align*}
	|q'''+2 \cdot f'''| = |2 \cdot f'''| &\leq 2 \cdot (2 \cdot f'')^{\frac{3}{2}} \\
	&= 4 \cdot \sqrt{2} \cdot (f'')^{\frac{3}{2}} \leq 4 \cdot \sqrt{2} \cdot (\underbrace{q''+f''}_{\geq 0} + f'')^{\frac{3}{2}}
\end{align*}
Thus,
\begin{align*}
	\left|\frac{1}{8}\cdot g'''\right|\leq 4 \cdot \sqrt{2} \cdot \left( \frac{1}{8} \cdot g''\right)^{\frac{3}{2}}\,.
\end{align*}
Pulling the factor $1/8$ out of the brackets, followed by multiplying both sides of the above inequality with $8$, we arrive at
\begin{align*}
	|g'''| \leq \underbrace{ 4 \cdot \sqrt{2} \cdot 8 \cdot \frac{1}{8^{\frac{3}{2}}} }_{\equiv 2} \cdot (g'')^{\frac{3}{2}}\,,
\end{align*}
which shows the proposition. q.e.d.

\largeparbreak

We summarize two simple further results that we will make use of.
\begin{cor}\label{cor:self}
The following hold for scalar or additive manipulations of self-concordant functions:
\begin{enumerate}
	\item Given a self-concordant function $g$. Then $c \cdot g$ is self-concordant for every parameter $c \in [1,\infty)$\,.
	\item Given two functions $g_1,g_2$, where $g_j$ is self-concordant in $\dom(g_j)$ $\forall j\in\lbrace 1,2\rbrace$. Then $g_1+g_2$ is self-concordant in $\dom(g_1) \cap \dom(g_2)$. If further $g_1$ or $g_2$ is strictly convex, then so is $g_1+g_2$.
\end{enumerate}
\end{cor}

The following corollary gives a quantitative result on strict convexity.
\begin{cor}\label{cor:strict_conv}
	Consider the self-concordant function
	$$ 	g:=f+\Gamma_{\Delta}\,. 	$$
	Let $f$ be a three times differentiable function, that is convex over $\Omega$. 
	
	Then, it holds
	$$ \lambda_{\text{min}}\big(\nabla^2 g(\bx)\big) \geq \frac{2}{\Delta^2} 	\quad\quad \forall \bx \in \Omega\,.$$
	Further, it holds
	$$ 	g(\bx) - \min\lbrace g \rbrace \leq 1/\Delta^2 \cdot \|\nabla g(\bx)\|_2^2\,. 	$$
\end{cor}
\noindent
\underline{Proof:} The first proposition follows from the analytic expression for the Hessian of $\GammDelta$.
For the second proposition we insert $\lambda_{\text{min}}(\nabla^2 g)\geq 1/\Delta^2$ into \eqref{eqn:BoundNewtonDecrement}. We then insert the obtained bound for the Newton decrement into \eqref{eqn:BoundOptimalityGap}. q.e.d.

\section{The method}
In this section we present our proposed method. It consists of three phases. The first phase applies Algorithm~\ref{algo:DampedNewton} to find the minimizer $\bx$ of a self-concordant function. The second and third phase apply Algorithm~\ref{algo:PathFollowNewton} subsequently to find minimizers of parametric self-concordant functions.

We now state the strictly convex self-concordant objective functions that the three phases minimize:
\begin{align*}
	\phi^{(1)}_{\phantom{\tau}}(\bx) &:= \GammD(\bx)\\
	\phi^{(2)}_\tau(\bx) &:= 8 \cdot \Big(\,\frac{2}{\tau} \cdot \hat{q}(\bx) + 2 \cdot \GammD(\bx)\,\Big)\\
	\phi^{(3)}_\pi(\bx) &:= \frac{8\cdot \tau_F}{\pi} \cdot \Big(\,\frac{2}{\tau_F} \cdot \hat{q}(\bx) + 2 \cdot \GammLR(\bx)\,\Big) + 16 \cdot \GammDelta(\bx)
\end{align*}
As we see, the third function is proportional to $\Phi(\bx)$ when $\pi = \pi_F$.
\largeparbreak

The section is organized as follows. We dedicate one subsection for each phase. In the subsections we prove that the objective function of the respective phase is actually self-concordant for the relevant range of parameters for $\tau$ and $\pi$. Then, we show how to construct a suitable initial guess. Certainly, except for the first phase, this is based on just taking the found solution from the former phase.

After having defined the initial guess, we analyze the number of iterations that Newton's method requires to converge from this initial guess to the global minimizer for the objective function of that phase. In the end, we provide a summary subsection with an overall complexity result, accuracy result, and a pseudo-code.

\subsection{Phase 1}
Certainly, $\phi^{(1)}$ is self-concordant since of log-barrier type. Importantly, we notice that $\nabla^2 \phi^{(1)}$ is diagonal. Hence, $\phi$ can be minimized by independently minimizing it in each component $x_j$, $j=1,...,n$\, of $\bx$. There result the following decoupled problems
\begin{equation}
	\tag{$\mathsf{BOX}_j$}\label{eqn:boxj}
	\begin{aligned}
		&\operatornamewithlimits{min}_{x \in (\ell_j,u_j)}& \quad \Gamma_j(x) :=&\quad-\log(x - x_{\mathsf{L},j}) - \log(x_{\mathsf{R},j}-x) \\
		& & &\quad-\log(\Delta + x) - \log(\Delta - x)
	\end{aligned}
\end{equation}
for each component $j=1,...,n$, where the limits of the open interval domain are
\begin{align*}
	\ell_j:= \max\lbrace\,-\Delta\,,\,x_{\mathsf{L},j}\rbrace\,, \quad\quad u_j:= \min\lbrace\,\Delta\,,\,x_{\mathsf{R},j}\rbrace\,.
\end{align*}
In fact, according to Wolfram Alpha, this scalar minimization problem even has a closed expression for the minimizer. However, this closed expression is complicated and it is unclear how to evaluate it in a numerically reliable way.

This is why in the following we consider the minimization of $\phi^{(1)}$ by separately minimizing the problems \eqref{eqn:boxj} for each $j=1,...,n$ by using Algorithm~\ref{algo:DampedNewton}.

\begin{theorem}\label{thm:Phase1IniGuess}
	Consider $\Gamma_j$ and its unique minimizer $x^\star \in (\ell_j,u_j) \subset \R$. Use the initial guess
	$$ x_0 := \frac{1}{2} \cdot (\ell_j + u_j)\,. $$
	Choose $\varepsilon>0$.
	
	Then, after $k \in \N$ iterations, Algorithm~\ref{algo:DampedNewton} finds an iterate $x_k$ that satisfies
	\begin{align*}
		\Gamma_j(x_k) - \Gamma_j(x^\star) \leq \varepsilon\,,
	\end{align*}
	where $k$ obeys to the bound
	$$	k \leq 64 + \log_2\big(1-\log_2(\varepsilon)\big). 	$$
\end{theorem}
\noindent
\underline{Proof:} Our proof works by constructing a worst-case problem instance that yields the largest amount of Newton iterations. We then simplify this instance without loss of generality.

Without loss of generality, we can choose $x_0 = 0$. This works by simply applying a shift to the ordinate. Further, without loss of generality, we can say that $u_j - \ell_j = 2 \cdot \Delta$. This is not a restriction because yet we have not imposed any further conditions on $\Delta,u_j,\ell_j$. After these simplifications, we arrive at the situation where
\begin{align*}
	\ell_j=-\Delta\,,\quad\quad u_j = \Delta\,,\quad\quad x_{\mathsf{L},j} \leq -\Delta\,,\quad\quad x_{\mathsf{R},j}\geq \Delta\,.
\end{align*}
Now we construct the worst case. Certainly, if $x_{\mathsf{L},j} = -\Delta$ and $x_{\mathsf{R},j} = \Delta$ then $x_0 = 0$ is already optimal. Hence, this would be the best case. Also, when $x_{\mathsf{L},j} \ll -\Delta$ and  $x_{\mathsf{R},j} \gg \Delta$, it appears that the shape of $\Gamma_j$ is dominated by $-\log(\Delta+x)-\log(\Delta+x)+const$, which again had a minimizer $x^\star\approx x_0$.

In conclusion, the worst case must appear when $x_{\mathsf{L},j}=-\Delta$ and $x_{\mathsf{R},j}\rightarrow \infty$. In this case, it holds
$$ \Gamma_j(x) = -2 \cdot \log(\Delta + x) - \log(\Delta-x) + c\,, $$
where $c \in \overline{\R}$ is a constant in $x$. From the optimality conditions we determine the minimizer $x^\star= \Delta/3$. By insertion, we find that $x^\star$ and our initial guess $x_0 = 0$ satisfy
\begin{align*}
	\Gamma_j(x^\star) &= c- 3 \cdot \log(\Delta)-\log\Big(\frac{32}{27}\Big) \,,\\
	\Gamma_j(x_0) &= c-3 \cdot \log(\Delta)\,.
\end{align*}
Using Theorem~\ref{thm:Newton_Minimization} shows the proposition:
$$ 	375 \cdot \big(\,\Gamma_j(x_0) - \Gamma_j(x^\star)\,\big)<64\,. 	$$
q.e.d.

\largeparbreak

We compute the minimizers $x_j$ of \eqref{eqn:boxj} for $j=1,2,...,n$ with tolerance
\begin{align}
	\varepsilon := \min\left\lbrace\,\left(\frac{\delta \cdot \Delta}{2048 \cdot \sqrt{n}}\right)^2\,,\,\frac{1}{36}\,\right\rbrace\,.	\label{eqn:DefVareps}
\end{align}
The obtained values $x_j$ are arranged into the vector $\bx_{[\mathsf{I}]} := (x_1,x_2,...,x_n)\t$, that forms the outcome of phase 1. This completes our work in phase 1. In the next subsection, we will continue with phase 2.

\largeparbreak

In preparation for the analysis of phase 2 and as a sequel to this subsection, in the lemma below we establish a bound for the gradient norm of the objective function at the numerically computed minimizer.

\begin{lem}\label{lem:BoundGradientPhase1}
Consider $\bx_{[\mathsf{I}]} \in \Omega$ as computed above, with $\varepsilon$ as defined in \eqref{eqn:DefVareps}.

Then, it holds
$$ 	\|\nabla \GammD(\bx_{[\mathsf{I}]})\|_2 \leq \Delta/64 \,.		$$
\end{lem}
\noindent
\underline{Proof:}
We start by observing
\begin{align*}
	\|\nabla\GammD(\bx_{[\mathsf{I}]})\|_2 \leq \sqrt{n} \cdot \|\nabla\GammD(\bx_{[\mathsf{I}]})\|_\infty = \sqrt{n} \cdot \operatornamewithlimits{max}_{1\leq j \leq n}\big\lbrace\, |\Gamma_j'(x_j)| \,\big\rbrace\,,
\end{align*}
where $x_j$ was the minimizer of \eqref{eqn:boxj}. In the following we analyze the solution $x_j$ and derive a bound for $|\Gamma_j'(x_j)|$.

We define
$$ 	x^\star_j := \operatornamewithlimits{argmin}_{\xi \in (\ell_j,u_j)} \big\lbrace\, \Gamma_j(\xi) \, \big\rbrace\,, 	$$
the exact minimizer of $\Gamma_j$.
From the worst-case analysis in the proof of Theorem~\ref{thm:Phase1IniGuess} we know that it holds $x^\star_j \in I_j$, with the interval
$$ 		I_j := \left[\, \ell_j + \frac{1}{3} \cdot (u_j - \ell_j) \ , \ u_j - \frac{1}{3} \cdot (u_j - \ell_j) \,\right]\,. $$
We also define a second, slightly wider, interval
$$ 		\hat{I}_j := \left[\, \ell_j + \frac{1}{4} \cdot (u_j - \ell_j) \ , \ u_j - \frac{1}{4} \cdot (u_j - \ell_j) \,\right]\,. 	$$
Clearly, $x^\star_j \in \hat{I}_j$. Notice that $\hat{I}_j$ is wider than $I_j$ in both directions by a margin of $1/12 \cdot (u_j-\ell_j)$.

We notice from insertion of critical points on $\hat{I}_j$ into the analytic expression of $\Gamma_j''$, that
$$ 		\frac{8}{(u_j - \ell_j)^2}\leq \Gamma_j''(\xi) \leq \frac{64}{(u_j - \ell_j)^2}\quad \quad \forall \xi \in \hat{I}_j\,. 	$$
Due to the strict convexity of $\Gamma_j$ and the optimality of $x^\star_j$, we can bound the objective from below by a parabola with second derivative $8/(u_j - \ell_j)^2$. Hence,
$$ 	\frac{1}{2} \cdot \frac{8}{(u_j - \ell_j)^2} \cdot |x_j - x^\star_j|^2 \leq \Gamma_j(x_j) - \Gamma_j(x^\star_j) \leq \varepsilon\,. 	$$
Taking the square-root and using $\varepsilon\leq 1/36$, we find
$$	|x_j - x^\star_j| \leq \sqrt{1/4} \cdot (u_j - \ell_j) \cdot \sqrt{\varepsilon} \leq 1/12 \cdot (u_j - \ell_j)\,. 	$$
Thus, $x_j \in \hat{I}_j$.

Using a Taylor series and inserting $\sqrt{\varepsilon}\leq \delta\cdot\Delta/2048$, we arrive at
\begin{align*}
		|\Gamma_j'(x_j)| 
		&\leq \underbrace{|\Gamma_j'(x^\star_j)|}_{=0} + \operatornamewithlimits{max}_{\xi \in \hat{I}_j}\big\lbrace\,|\Gamma_j''(\xi)|\,\big\rbrace \cdot |x_j - x^\star_j| \\
		&\leq \frac{64}{(u_j - \ell_j)^2} \cdot |x_j -x^\star_j| \leq \frac{32}{u_j - \ell_j} \cdot \underbrace{2 \cdot \frac{|x_j - x^\star_j|}{u_j-\ell_j}}_{\leq \sqrt{\varepsilon}}\\
		& \leq \frac{32}{u_j-\ell_j} \cdot \frac{\delta \cdot \Delta}{2048 \cdot \sqrt{n}} \leq \frac{1}{\sqrt{n}} \cdot \frac{\Delta}{64}\,.
\end{align*}
q.e.d.

\subsection{Phase 2}
After $\bx_{[\mathsf{I}]}$ has been computed, we start with phase 2. In this phase we compute the minimizer of $\phi^{(2)}_{\tau_F}$\,.

We start by showing the self-concordance of $\phi^{(2)}_\tau$.

\begin{cor}
	The parametric function $\phi^{(2)}_\tau$ is self-concordant in $\bx$ over $\Omega$ for every parameter $\tau \in [\tau_F,\infty)$.
\end{cor}
\noindent\underline{Proof:} In Section~1 we required that the function
$$	\psi = \hat{q} + \tau_F/2 \cdot \GammLR 	$$
is convex over $\Omega$. Hence, also 
$$ 	\tilde{\psi}_\tau := 2/\tau \cdot \hat{q} + \GammLR + \GammDelta \equiv 2/\tau \cdot \hat{q} + \GammD	$$
is convex over $\Omega$ for every parameter $\tau \in [\tau_F,\infty)$. Applying Lemma~\ref{lem:self_from_conv} to the function $\tilde{\psi}_\tau$ for a respective value or $\tau$, we find that
$$ 	8 \cdot \Big(\,2/\tau \cdot \hat{q} + 2 \cdot \GammD \,\Big)\equiv \phi^{(2)}_\tau 	$$
is self-concordant over $\Omega$ for every parameter $\tau \in [\tau_F,\infty)$ as well. q.e.d.

\largeparbreak

We want to use Algorithm~\ref{algo:PathFollowNewton} to minimize $\phi^{(2)}_\tau$ for decreasing values of $\tau$. To start the iteration, we want to use the initial guess $\bx_0 := \bx_{[\mathsf{I}]}$. But, we must find a suitable value for $\tau_0$ such that the following requirement from Theorem~\ref{thm:NewtonPathFollow} is satisfied:
\begin{align}
	\phi^{(2)}_{\tau_0}(\bx_0) - \min\lbrace\,\phi^{(2)}_{\tau_0}\,\rbrace \leq 1/4\,. \label{eqn:BoundConditionPhase2Iniguess}
\end{align}

Certainly, there is hope that there is a suitably large value for $\tau_0$ to achieve this bound, because
\begin{align*}
	\operatornamewithlimits{limes}_{\tau\rightarrow\infty}\quad \phi^{(2)}_\tau = 16 \cdot \phi^{(1)}\,.
\end{align*}
Hence, the optimality of $\bx_{[\textsf{I}]}$ implies that it is also almost optimal for $\phi^{(2)}_{\tau_0}$, subject that $\tau_0$ is very large. Below, we derive such a suitably large value for $\tau_0$: 

\largeparbreak

We apply Corollary~\ref{cor:strict_conv}, that gives the bound
\begin{align}
\phi^{(2)}_{\tau_0}(\bx_0) - \min\lbrace\,\phi^{(2)}_{\tau_0}\,\rbrace \leq \frac{1}{\Delta^2} \cdot \|\nabla\phi^{(2)}_{\tau_0}(\bx_0)\|_2^2\,. \label{eqn:BoundPhi2Ini}
\end{align}
Inserting \eqref{eqn:BoundPhi2Ini} into the requirement \eqref{eqn:BoundConditionPhase2Iniguess} and taking the square, we find that it must hold:
$$ 	\frac{1}{\Delta} \cdot \|\nabla\phi^{(2)}_{\tau_0}(\bx_0)\|_2 \leq \frac{1}{2} 	$$ 
Remembering $\bx_0 = \bx_{[\mathsf{I}]}$ and using the bound from Lemma~\ref{lem:BoundGradientPhase1}, we can bound the norm of the gradient as
\begin{align*}
	\|\nabla\phi^{(2)}_{\tau_0}(\bx_0)\|_2
	&\leq \Big\|8 \cdot \Big(\,2/\tau_0 \cdot \nabla\hat{q}(\bx^\star)+2 \cdot \underbrace{\nabla \GammD(\bx_{[\mathsf{I}]})}_{\leq \Delta/64}\Big)\Big\|_2\\
	&\leq \frac{16}{\tau_0} \cdot \Big(\|\bQ\|_2 \cdot (\|\bxL\|_2+\|\bxR\|_2)+\|\bc\|_2\Big) + \Delta/4 \,.
\end{align*}
Hence, condition \eqref{eqn:BoundConditionPhase2Iniguess} is satisfied when we bound the first term above by $\Delta/4$. We can achieve that by defining
\begin{align*}
	\tau_0 := \frac{64}{\Delta} \cdot \Big(\|\bQ\|_2 \cdot (\|\bxL\|_2+\|\bxR\|_2)+\|\bc\|_2 \Big) \,.
\end{align*}

In summary, we apply Algorithm~\ref{algo:PathFollowNewton} to follow the path of minimizers of $\phi^{(2)}_\tau$. $\tau$ starts from $\tau_0$ with initial guess $\bx_0 := \bx^\star$, and $\tau$ decreases until $\tau=\tau_F$. The obtained minimizer for $\phi^{(2)}_{\tau_F}$ from phase 2 we call $\bx_{[\mathsf{II}]}$\,. For the computation of $\bx_{[\mathsf{II}]}$, the tolerance $\varepsilon = 1/4$ is sufficiently accurate.

\subsection{Phase 3}
After $\bx_{[\mathsf{II}]}$ has been computed, we start phase 3. In this phase we compute the minimizer of $\phi^{(3)}_{\pi_F}$.

We start again by showing self-concordance for $\phi^{(3)}_\pi$. Remember that in Section~1 we stated the requirement $\tau_F\geq \pi_F$.

\begin{cor}
	The parametric function $\phi^{(3)}_\pi$ is self-concordant in $\bx$ over $\Omega$ for every parameter $\pi \in (0,\tau_F]$.
\end{cor}
\noindent\underline{Proof:} In Section~1 we required that the function
$$	\psi = \hat{q} + \tau_F/2 \cdot \GammLR 	$$
is convex over $\Omega$. Hence, also 
$$ 	\tilde{\psi} := 2/\tau_F \cdot \hat{q} + \GammLR	$$
is convex over $\Omega$. Applying Lemma~\ref{lem:self_from_conv} to the function $\tilde{\psi}$, we find that
$$ 	8 \cdot \Big(\,2/\tau_F \cdot \hat{q} + 2 \cdot \GammLR \,\Big) 	$$
is self-concordant over $\Omega$. From Corollary~\ref{cor:self} we know that scaling of a function with constant $\geq 1$ does not destroy self-concordance. Hence, the function
$$  8 \cdot \frac{\tau_F}{\pi} \cdot \Big(\,2/\tau_F \cdot \hat{q} + 2 \cdot \GammLR \,\Big) 	$$
is self-concordant over $\Omega$ for every parameter $\pi \in (0,\tau_F]$.
Finally, we add the self-concordant function $16 \cdot \GammDelta$ to it, which, according to Corollary~\ref{cor:self}, yields again a self-concordant function:
$$ 	8 \cdot \frac{\tau_F}{\pi} \cdot \Big(\,2/\tau_F \cdot \hat{q} + 2 \cdot \GammLR \,\Big) + 16 \cdot \GammDelta \equiv \phi^{(3)}_\pi 	$$
q.e.d.

\largeparbreak

We find that for $\pi = \tau_F$ it holds $\phi^{(3)}_\pi \equiv \phi^{(2)}_{\tau_F}$. Hence, the minimizer $\bx_{[\mathsf{II}]}$ is the minimizer of $\phi^{(3)}_\pi$ for $\pi = \pi_0 := \tau_F$. Hence, we can apply Algorithm~\ref{algo:PathFollowNewton} to follow the path of minimizers for $\phi^{(3)}_\pi$. The iteration starts from $\pi = \pi_0$ with initial guess $\bx_0 := \bx_{[\mathsf{II}]}$. $\tau$ decreases iteratively until $\pi = \pi_F$ is reached.

The returned solution of Algorithm~\ref{algo:PathFollowNewton} we call $\bx_{[\mathsf{III}]}$. This vector is a numerical solution for the minimizer of $\phi^{(3)}_{\pi_F}$. It satisfies the accuracy condition
$$ 	\phi^{(3)}_{\pi_F}(\bx_{[\textsf{III}]}) - \min\lbrace\,\phi^{(3)}_{\pi_F}\,\rbrace \leq \varepsilon\,, 	$$
when calling Algorithm~\ref{algo:PathFollowNewton} with the tolerance parameter $\varepsilon>0$. Eventually, we probably rather want a solution that satisfies
$$ 	\Phi(\bx_{[\textsf{III}]}) - \min\lbrace\,\Phi\,\rbrace \leq \tol 	$$
for the objective $\Phi$ of \eqref{eqn:NCQPTR} and a user-demanded tolerance $\tol>0$\,.
Since $\Phi = 16/\pi_F \cdot \phi^{(3)}_{\pi_E}$, we achieve the demanded tolerance by choosing $\varepsilon := \tol \cdot \pi_F/16$\,.

After $\bx_{[\mathsf{III}]}$ has been computed, it is returned as the solution of \eqref{eqn:NCQPTR}, and our proposed method terminates with success.

\subsection{Algorithm and complexity in summary}

\paragraph{The entire algorithm} is stated in Algorithm~\ref{algo:Proposed}. The definitions of the quadratic function $\hat{q}$ and the logarithmic barrier-functions $\GammLR,\GammDelta,\GammD$ are given in \eqref{eqn:hatq}---\eqref{eqn:hatGamma}.

\begin{algorithm}
	\caption{Proposed method}
	\label{algo:Proposed}
	\begin{algorithmic}[1]
		\Procedure{Solver}{$\bQ,\bc,\bxL,\bxR,\Delta,\tau_F,\pi_F,\tol$}\vspace{2mm}
		\State \textit{// - - - Phase 1 - - -}
		\State $\delta := \operatornamewithlimits{min}_{1 \leq j \leq n}\Big\lbrace \ \min\lbrace x_{\textsf{R},j},\Delta\rbrace-\max\lbrace-\Delta,x_{\textsf{L},j}\rbrace \ \Big\rbrace$
		\State $\varepsilon := \min\Big\lbrace\,\big(\frac{\delta \cdot \Delta}{2048 \cdot \sqrt{n}}\big)^2\,,\,1/36\,\Big\rbrace$
		\For{$j=1,...,n$}
		\State \textit{// solve problem \eqref{eqn:boxj}}
		\State $x_0 := 1/2 \cdot (\max\lbrace-\Delta,x_{\mathsf{L},j}\rbrace + \min\lbrace\Delta,x_{\mathsf{R},j}\rbrace)$
		\State $g(x) := -\log(x-x_{\mathsf{L},j})-\log(x_{\mathsf{R},j}-x)-\log(\Delta+x)-\log(\Delta-x)$
		\State $x_j := $\textsc{DampedNewton}($g,x_0,\varepsilon$)
		\EndFor
		\State $\bx_{[\textsf{I}]} := (x_1,x_2,...,x_n)\t$\vspace{2mm}
		\State \textit{// - - - Phase 2 - - -}
		\State $f(\bx) := 16 \cdot \hat{q}(\bx)$
		\State $\Gamma(\bx) := 16 \cdot \GammD(\bx)$
		\State $\tau_0 := 64/\Delta \cdot \Big(\|\bQ\|_2 \cdot (\|\bxL\|_2+\|\bxR\|_2)+\|\bc\|_2\Big)$
		\State $\bx_{[\mathsf{II}]}:=$\textsc{PathFollowNewton}($f,\Gamma,\bx_{[\textsf{I}]},\tau_0,\tau_F,1/4$)\vspace{2mm}
		\State \textit{// - - - Phase 3 - - -}
		\State $f(\bx) := 16 \cdot \big(\hat{q}(\bx)+\tau_F \cdot \GammLR(\bx)\big)$
		\State $\Gamma(\bx) := 16 \cdot \GammDelta(\bx)$
		\State $\pi_0 := \tau_F$\,,\quad\quad$\varepsilon:=\tol \cdot \pi_F / 16$
		\State $\bx_{[\textsf{III}]} := $\textsc{PathFollowNewton}($f,\Gamma,\bx_{[\mathsf{II}]},\pi_0,\pi_F,\varepsilon$)\vspace{2mm}
		\State $\hbx := \bx_{[\textsf{III}]}$
		\State \Return $\hbx$
		\EndProcedure
	\end{algorithmic}
\end{algorithm}

\paragraph{The complete complexity result}
We now give the complexity of the method. For a succinct presentation we define the problem size
\begin{align*}
L := 	&\quad1+\log(1+\|\bQ\|_2) + \log(1+\|\bc\|_2) \\
		&\quad+\log(1+\|\bxL\|_2)+ \log(1+\|\bxR\|_2) + \log(1+\Delta+1/\Delta) \\
		&\quad + \log(n) - \log(\pi_F)+\log\big(1-\log(\tol)\big)+\log\big(1-\log(\delta)\big)\,.
\end{align*}

The first phase has complexity $\cO(n \cdot L)$. This is because the solution of each subproblem \eqref{eqn:boxj} is a one-dimensional problem that converges in $\cO(L)$ Newton iterations. But, there are no linear systems that must be solved. Hence, this phase can actually be neglected in complexity terms because it is dominated by the cost of the subsequent phases.

The second phase uses a barrier-function $\Gamma$ with $|\Gamma| \in \cO(n)$. Using Theorem~\ref{thm:NewtonPathFollow}, we can bound the total number of Newton iterations (outer times inner iterations). There total number of iterations is bounded by the complexity
$$ 	\cO(\sqrt{n} \cdot L) 	\,.$$
The same complexity bound as for the second phase holds also for the third phase.

\largeparbreak
We now summarize the overall complexity: If we assume that the positive definite linear equation systems of size $n \times n$ in Algorithm~\ref{algo:DampedNewton} line~3 are solved via dense Cholesky factorization in $\cO(n^3)$ respectively, then the total time-complexity of Algorithm~\ref{algo:Proposed} is
$$ 	\cO(n^{3.5} \cdot L) 	\,.$$

\paragraph{Accuracy of the solution}
The returned solution vector $\hbx \in \Omega$ of Algorithm~\ref{algo:Proposed} satisfies
\begin{align*}
	\Phi(\hbx) - \operatornamewithlimits{min}_{\bx \in \Omega}\big\lbrace\,\Phi(\bx)\,\big\rbrace \leq \tol\,.
\end{align*}
This result has been derived in the end of Section~3.3\,.

\section{Outlook}
We have introduced a particular unconstrained minimization problem \eqref{eqn:NCQPTR} with a non-convex quadratic function and barrier-functions, that account for box-constraints and trust-region constraints. We motivated that solutions to such problems can provide efficient search-directions for nonlinear programming methods.

We then presented a weakly polynomial efficient algorithm for the solution of \eqref{eqn:NCQPTR}. This method is based on three phases, of which each applies a Newton technique for the minimization of self-concordant functions. The method is a primal interior-point method.

\largeparbreak

While theoretically efficient, our method uses a short-step path-following mechanism. This mechanism makes a defensive choice for the geometric reduction factor $\sigma <1$. This is known to result in very many linear equation systems that must be solved, which is computationally prohibitive and hence yet impractical.

Practical primal methods simply choose $\sigma$ much smaller. In this case, polynomial complexity results for the worst case are still available, however they are much worse. In \cite{Boyd} the authors suggest values in the range $1/\sigma \in [10,20]$ as reasonable practical choices. They present numerical experiments for linear programming, where this strategy results on average in $40$ to $60$ Newton iterations in total. This constitutes a very desirable practical complexity; however, it is an open question whether such good performance will also appear in our case, where we solve functions with an additional non-convex quadratic term.

\largeparbreak

Experience suggests that primal-dual methods may be faster in practice than purely primal methods. The method would immediately become practical if the short-step primal path-following could be replaced by a primal-dual path-following with a heuristic like in Mehrotra's predictor-corrector method \cite{Mehrotra}. This heuristic has proven to be efficient when the primal and dual KKT residual are linear functions.

But in contrast, in this paper we consider functions where, depending on the formulation, the primal and dual KKT residuals (or in alternative an augmented complementarity residual) are nonlinear. For this reason, further research and trial is necessary in order to find a practically suitable path-following scheme.

\FloatBarrier

\bibliography{NCQP_TR_bib}
\bibliographystyle{plain}

\end{document}